\documentclass[12pt,reqno]{amsart} 
\usepackage{amsthm,amsmath,amssymb,amsfonts,amscd,latexsym,verbatim,xypic}
\newtheorem{Theorem}{Theorem}[section]
\newtheorem{Lemma}[Theorem]{Lemma}
\newtheorem{Proposition}[Theorem]{Proposition}

\newtheorem{Remark}[Theorem]{Remark}

\newtheorem{Question}[Theorem]{Question}
\newcommand{\demo}{\noindent {\sc Proof.}\;}
\newcommand{\be} {\begin{equation}}
\newcommand{\ee} {\end{equation}}
\newcommand{\bea} {\begin{eqnarray}}
\newcommand{\eea} {\end{eqnarray}}
\newcommand{\lp}  {\left(}
\newcommand{\rp}  {\right)}
\newcommand{\eqdef} {\stackrel{\rm def}{=}}
 
\newcommand{\ZZ}{\mathbb{Z}}

\newcommand{\CC}{\mathbb{C}}
\newcommand{\PP}{\mathbb{P}}

\newcommand{\caR}{\mathcal{R}}
\newcommand{\bR}{\bar{\caR}}
\newcommand{\cB}{\mathcal{B}}
\newcommand{\fS}{\mathfrak{S}}

\newcommand{\ux}{\mathbf x}
\newcommand{\al}{\alpha} 
\newcommand{\ph}{\phi} 
\newcommand{\la}{\lambda} 
\newcommand{\La}{\Lambda} 
\newcommand{\de}{\delta}
\newcommand{\De}{\Delta}
\newcommand{\si}{\sigma}   
\newcommand{\ep}{\epsilon} 
\newcommand{\ga}{\gamma} 
\newcommand{\om}{\omega}
\begin{document} 
\title[On the invariants of the binary quintic]
{A computational solution to a question by
Beauville on the invariants of the binary quintic}
\author[Abdesselam]
{Abdelmalek Abdesselam}
\maketitle

\parbox{12cm}{\small 
{\sc Abstract.} 
We obtain an alternate proof of an injectivity result by Beauville
for a map from the moduli space of quartic del Pezzo
surfaces to the set of conjugacy classes of certain subgroups
of the Cremona group. This amounts to showing that a projective
configuration of five distinct unordered points on the line
can be reconstructed from its five projective four-point
subconfigurations. This is done by reduction to a question
in the classical invariant theory of the binary quintic, which
is solved by computer-assisted methods. More precisely,
we show that six specific invariants of degree $24$,
the construction of which was explained to us by Beauville,
generate all invariants the degrees of which are divisible by $48$. 
}

\vspace{5mm} 

\mbox{\small AMS subject classification (2000): 14-04; 68W30; 12Y05;
14E07; 20G05}

\medskip 

\parbox{12cm} 


\section{Introduction}

Throughout this article, our base field will be $\CC$.
Let ${\rm Cr}$ denote the Cremona group of birational
transformations of $\PP^2$. To any element $S$ in the moduli
space of quartic del Pezzo surfaces, one can naturally associate
an element $G_S$ in the set of conjugacy classes of subgroups
isomorphic to $(\ZZ/2)^4$ inside ${\rm Cr}$.
This construction was considered in the recent 
work of Beauville~\cite{Beauville2}; among the results he proves
therein, one finds the following statement (loc. cit., Prop. 4.2).

\begin{Proposition}\label{Bmapinj}
The map $S\rightarrow G_S$ is injective.
\end{Proposition}

In the mentioned article, this result was obtained
by an elegant geometric argument,
using an idea of Iskovskikh~\cite{Iskovskikh}.
However, in an earlier version of the same work~\cite{Beauville1},
the weaker statement of generic injectivity was obtained by a
radically different approach, with a flavor of classical
invariant theory.
The purpose of the present article is to push this second approach
to completion, and show that it leads to a strengthening
of Proposition \ref{Bmapinj}, which is Theorem \ref{Reconst} below.
Our proof however is a {\em computer-assisted} one, since it relies on
rather heavy calculations using the Maple software.

In this classical invariant theoretic setting, the quite pretty
`reconstruction problem' that needs to be solved is the following.
Let $\La=\{\la_1,\ldots,\la_5\}$ be a set
of five distinct unordered points on $\PP^1$, and consider
the quintic
\be
\caR\eqdef\prod\limits_{1\le l\le 5} (z-j_l)
\label{jquintic}
\ee
where $j_l$ is the well-known $j$-invariant
of the four-point subset 
\[
\{\la_1,\ldots,\widehat{\la_l},
\ldots,\la_5\}\ \ .
\] 

\begin{Question}
\label{Qreconst}
Does the quintic $\caR$ uniquely determine
the $SL_2$ orbit of $\La$?
In other words, can one reconstruct the projective
configuration of a five-point set on the line from
the projective configurations of its four-point subsets?
\end{Question}

\begin{Remark}
Most of the difficulty here stems from the lack of any
ordering information, as well as the possibility of deforming
each of the five four-point pictures by a priori {\em unrelated}
homographies.
\end{Remark}

Theorem \ref{Reconst} below gives
an affirmative answer to this question, and
also implies Proposition \ref{Bmapinj}.
Indeed, by considering a homogenized version of the quintic $\caR$,
one is naturally led to the construction of six invariants
$\cB_0,\ldots,\cB_5$ of the binary quintic corresponding to the quintuple
$\La$, all of degree $24$.
Question \ref{Qreconst} is then solved by reduction to the following one.

\begin{Question}
\label{Qinv}
Is there a strictly positive integer $d_0$, divisible by $24$,
such that for all multiples $d$ of $d_0$, all invariants of
the binary quintic which have
degree $d$ can be polynomially expressed in terms of
the invariants
$\cB_0,\ldots,\cB_5$?
\end{Question}

Note that there are $7$ linearly independent invariants
in degree $24$, therefore if such a $d_0$ exists it has to be no smaller
than $48$. Theorem \ref{Thm48} below, shows that $d_0=48$ indeed does the job,
providing a positive answer to Question \ref{Qinv}.
Trying to understand the intriguing rather high degree at which this
phenomenon first occurs was our primary motivation
for the present work.

More precise statements of our results as well as
the detailed explanation of the steps in our calculations
will be given in Section \ref{Invthcal};
after the necessary material from the classical invariant theory
of binary forms is recalled in Section \ref{Prelim}.
In Section \ref{Geom}, we will briefly relate our results with Beauville's.
Finally, Section \ref{Todo} will outline some suggestions for further work.

\section{Preliminaries on the classical invariant theory of binary forms} 
\label{Prelim}

\subsection{Covariants, invariants and symmetric functions of root differences}
The following material is classical. However, it is recalled here
firstly for the convenience of the reader, and secondly in order
to fix the numerical normalization of the invariants we will be considering.

A binary form of order $p$ is a homogeneous polynomial
\be
F(\ux)=\sum\limits_{i=0}^p
a_i x_1^{p-i} x_2^i
\label{binform}
\ee
of degree $p$ in the variables $\ux\eqdef(x_1,x_2)$.
A matrix
\[
g=\lp
\begin{array}{ll}
g_{11} & g_{12}\\
g_{21} & g_{22}
\end{array}
\rp
\]
in $GL_2$ acting on the variables by
\[
\ux\rightarrow g\ux=(g_{11}x_1+g_{12}x_2, g_{21}x_1+g_{22}x_2)\ \ ,
\]
induces a transformation $F\rightarrow gF$
on the coefficients of the binary form $F$, by forcing the
equality
\[
(gF)(\ux)\eqdef F(g^{-1}\ux)\ .
\]
A covariant of $F$, of degree $d$, order $r$ and weight
$\om$, is a polynomial $C(F,\ux)=C(a_0,\ldots,a_p;x_1,x_2)$,
homogeneous of total degree $d$ in $a_0,\ldots,a_p$,
and homogeneous of total degree $r$ in $x_1,x_2$,
such that for any $g$ in $GL_2$,
\be
C(gF, g\ux)=({\rm det}\ g)^{-\om} C(F,\ux)\ .
\ee
One has a simple relation between $p$, $r$, $d$, and $\om$:
\be
d p=2\om+r\ .
\label{Degweight}
\ee
An invariant $I=I(a_0,\ldots,a_p)=I(F)$ simply is a
covariant of order zero.
For two pairs of variables $b=(b_1,b_2)$ and $c=(c_1,c_2)$
which can be thought of as the homogeneous coordinates
of two generic points in $\PP^1$, following the elegant classical
notation, we write
\[
(bc)\eqdef
\left|
\begin{array}{ll}
b_1 & c_1 \\
b_2 & c_2
\end{array}
\right|\ \ .
\]
In terms of its homogeneous roots
$\xi_1,\ldots, \xi_p$, the form
$F$ can therefore be written as
\[
F(\ux)=(\ux\xi_1)\ldots(\ux\xi_p)\ .
\]
Now one has the following classical result
(see~\cite[p. 97]{Elliott} or~\cite{KungR}).

\begin{Proposition}
\label{Rootprop}

\noindent 1)
Every symmetric polynomial in the pairs of variables $\xi_1,\ldots,\xi_p$
which is a linear combination of expressions of the form
\[
\prod\limits_{{1\le i,j\le p}\atop{i\neq j}}
(\xi_i \xi_j)^{k_{ij}}
\ \times\ \prod\limits_{1\le i\le p} (\ux\xi_i)^{l_i}
\]
where the $k$'s and the $l$'s are nonnegative integers satisfying
\[
\sum\limits_{1\le i\le p} l_i=r\ \ ,
\]
\[
\forall i, \ \ \sum\limits_{{1\le j\le p}\atop{j\neq i}}
(k_{ij}+k_{ji})\ +\ l_i=d\ \ ,
\]
and
\[
\sum\limits_{{1\le i,j\le p}\atop{i\neq j}} k_{ij}=\om\ \ ,
\]
is an (irrational) expression for a covariant of $F$, of degree $d$,
order $r$ and weight $\om$.

\noindent 2)
Conversely any covariant $C$ of $F$ can be so written.
\end{Proposition}

Note that the proposition has an obvious generalization to the
case of simultaneous covariants of more than one form.
For example, if one considers two binary forms
\[
F(\ux)=\sum\limits_{i=0}^p
a_i x_1^{p-i} x_2^i
=(\ux\xi_1)\ldots(\ux\xi_p)\ \ ,
\]
and
\[
G(\ux)=\sum\limits_{i=0}^q
b_i x_1^{q-i} x_2^i
=(\ux\eta_1)\ldots(\ux\eta_q)\ \ ,
\]
the resultant, which is a joint invariant of $F$ and $G$, is
\bea
{\rm Res}(F,G) & \eqdef & \prod\limits_{{1\le i\le p}\atop{1\le j\le q}}
(\xi_i \eta_j)\\
 & = &
\left|
\begin{array}{cccccc}
a_0 & \ldots & \ldots & a_p & & 0 \\
 & \ddots & & & \ddots & \\
0 & & a_0 & \ldots & \ldots & a_p \\
b_0 & \ldots & \ldots & b_q & & 0 \\
 & \ddots & & & \ddots & \\
0 & & b_0 & \ldots & \ldots & b_q \\
\end{array}
\right|
\label{Sylvres}
\eea
the usual Sylvester $(p+q)\times(p+q)$-determinant formula.
Likewise the discriminant of a form $F$ is by definition
the invariant
\bea
{\rm Disc}(F) & \eqdef & \prod\limits_{1\le i<j\le p} (\xi_i \xi_j)^2\\
 & = & \frac{(-1)^{\frac{p(p-1)}{2}}}{p^{p-2}}
{\rm Res}\lp
\frac{\partial F}{\partial x_1},\frac{\partial F}{\partial x_2}
\rp\ \ .
\eea

We now need to recall the classical notion of transvectant
(or the ``Uebereinanderschiebung'' of~\cite[\S 1]{Gordan1}), which
allows the formation of new covariants from old ones,
and the formulation of quick yet {\em precise} definitions for
those used in Section \ref{Invthcal}.
If $F$ is a binary form of order $p$ and $G$ a binary form of order $q$,
the $k$-th transvectant of $F$ and $G$ is
\be
(F,G)_k = \frac{(p-k)!(q-k)!}{p!\ q!}  \, 
\sum\limits_{i=0}^k \, (-1)^i \binom{k}{i} \, 
\frac{\partial^{\,k} F}{\partial x_1^{k-i} \, \partial x_2^i } \; 
\frac{\partial^{\,k} G}{\partial x_1^i \, \partial x_2^{k-i}} 
\; .
\ee

\subsection{Invariants of the binary quartic}
The ring of invariants of a generic binary form $F$ of
order $p$ as in (\ref{binform})
is denoted by $\CC[a_0,\ldots,a_p]^{SL_2}$ or simply
$\CC[F]^{SL_2}$. It is given the grading by the degree
in the coefficients of $F$. The graded component of degree
$d$ is denoted by $\CC[F]_d^{SL_2}$.
For a binary quartic, more conveniently written
\be
Q(\ux)= q_0 x_1^4 + 4 q_1 x_1^3 x_2 + 6 q_2 x_1^2 x_2^2
+ 4 q_3 x_1 x_2^3 + q_4 x_2^4\ \ ,
\label{quartic}
\ee
the ring of invariants has been know
since the time of Boole and Cayley.
Following~\cite[p. 189]{Salmon} it can be described as
\be
\CC[Q]^{SL_2}=\CC[S,T]
\ee
where
\bea
S(Q) & \eqdef & \frac{1}{2}(Q,Q)_4 \\
 & = & q_0 q_4-4 q_1 q_3 +3 q_2^2
\label{Sformula}
\eea
is of degree $2$ and weight $4$, and
\bea
T(Q) & \eqdef & \frac{1}{6}\lp Q,(Q,Q)_2\rp_4 \\
 & = & q_0 q_2 q_4+ 2 q_1 q_2 q_3 - q_2^3
- q_0 q_3^2 - q_1^2 q_4
\label{Tformula}
\eea
is of degree $3$ and weight $6$;
besides, $S$ and $T$ are algebraically independent.
One also has the weight $12$ invariant
\be
{\rm Disc}(Q)= 2^8 \lp
S(Q)^3-27 T(Q)^2
\rp
\ee
as is readily checked on the canonical form written
with obvious notation
\be
Q(\ux)=(\ux 0)(\ux 1)(\ux \infty)(\ux \la)\ \ .
\ee
The classical $j$-invariant of the four-point set in $\PP^1$
corresponding to the roots of $Q$ is
\be
j(Q)\eqdef
\frac{S(Q)^3}{S(Q)^3-27 T(Q)^2}
=\frac{4}{27}\times
\frac{(\la^2-\la+1)^3}{\la^2 (\la-1)^2}\ \ .
\label{jinv}
\ee

\subsection{Invariants of the binary quintic}
For a binary quintic $F$, the description of the ring of invariants 
was completed by Hermite~\cite{Hermite1} building on
the previous work of Cayley and Sylvester.
Again according to~\cite[pp. 227--234]{Salmon},
on can describe it as
\be
\CC[F]^{SL_2}=\CC[J,K,L,H]/{\rm Relation}
\ee
where the invariants $J$, $K$, $L$, and $H$ are respectively
of degree 4, 8, 12, and 18; and there is a unique relation
between them in degree $36$ expressing $H^2$ in terms
of $J$, $K$, and $L$.
More precisely, one can make the following choices for the
generators. First, define the covariants
\bea
C_1 & \eqdef & (F,F)_4 \ ,\\
C_2 & \eqdef & (F, C_1)_2 \ ,\\
C_3 & \eqdef & (C_2,C_2)_2 \ ,\\
C_4 & \eqdef & (C_2,C_1)_2 \ .
\eea
Now the invariants are defined as
\bea
J & \eqdef & -\frac{1}{2} (C_1,C_1)_2 \ ,\label{Jdef}\\
K & \eqdef & \frac{1}{8} (C_1,C_3)_2 \ ,\\
L & \eqdef & \frac{1}{96} (C_3,C_3)_2 \ ,\\
H & \eqdef & -\frac{1}{384}
\lp (C_4,C_3)_1,(C_1,C_4)_1 \rp_1 \ .\label{Hdef}
\eea

Note that the full-fledged 
{\em Cartesian expressions}
for these invariants as linear combinations of monomials in the coefficients
of $F$ are quite complicated. Indeed, $J$, $K$, $L$, $H$ respectively have
$12$, $68$, $228$, and $848$ terms.
In order to calculate with invariants of the quintic $F$,
we will sometimes find it convenient to use the Sylvester
canonical form
\be
F(\ux)=u x_1^5 + v x_2^5- w (x_1+x_2)^5\ .
\label{Sylvcan}
\ee
Indeed, every form $F$ in the affine open set $\{L\neq 0\}$
can be written as the sum of the fifth powers of three
nonproportional linear forms.
By taking these points in the dual $\PP^1$ to
$0$, $1$, and $\infty$, one sees that such
an $F$ is in the $SL_2$ orbit of a form as in (\ref{Sylvcan}).
The reason for this is that
\be
L=-\frac{1}{2^4.3^5} {\rm Disc}({\rm Can}(F))
\ee
where the canonizant of $F$ is
\be
{\rm Can}(F)\eqdef - C_2 \ ,
\ee
or, in classical symbolic notation~\cite[\S 2]{AC2},
\be
{\rm Can}(F)= (ab)^2(ac)^2(bc)^2 a_\ux b_\ux c_\ux\ .
\ee
The above linear forms correspond to the distinct linear factors
of ${\rm Can}(F)$ (see e.g.~\cite[pp. 153--156]{Salmon} or~\cite{Kung}).
The point of this discussion is that any identity in the ring
$\CC[F]^{SL_2}$ can be checked by specialization to this canonical form.
The fundamental invariants will then be given by the remarkably simple
expressions:
\bea
J & = & (uv+uw+vw)^2-4uvw(u+v+w) \ ,\label{JSylv}\\
K & = & u^2 v^2 w^2(uv+uw+vw) \ ,\\
L & = & u^4 v^4 w^4 \ ,\\
H & = & u^5 v^5 w^5 (u-v)(u-w)(v-w) \ .\label{HSylv}
\eea
\begin{Remark}
The latter explain our choice of numerical normalization
in (\ref{Jdef}--\ref{Hdef}).
The explanation of the construction scheme we used
based on the covariants $C_1,\ldots,C_4$, is that it is
the most straightforward way to build, as a `Lego game', the
`Feynman diagrammatic' expression of the four invariants
(see~\cite[\S 6]{AC1} and~\cite[p. 120]{Kempe}).
The sums over `Wick contractions' involved in each of the
transvectant
operations produce, up to symmetry, only one graph.
Also note that (\ref{JSylv}--\ref{HSylv}) exactly agree with Salmon's
conventions~\cite{Salmon}, except for
the Hermite invariant $H$ which differs in sign and notation.
The invariants given by Gordan~\cite[\S 9]{Gordan1}
are different from the ones we
used here.
\end{Remark}

The unique relation, which can easily be checked
using (\ref{JSylv}--\ref{HSylv}),
is
\[
16H^2=
-432 L^3
-72 L^2KJ
+8 LK^3
-2 LK^2J^2
+L^2J^3
+K^4J\ .
\]
\be
{\ }
\label{Relation}
\ee

The dimension of a graded component of degree $d$ which is divisible
by $4$ can easily be calculated by solving for the nonnegative solutions
of an elementary diophantine equation. Indeed, because of the relation
(\ref{Relation}) one simply has to count the monomials in the algebraically
independent invariants $J$, $K$, $L$, with the given degree.
In sum,
\be
{\rm dim}\lp
\CC[F]_d^{SL_2}\rp
=\nu(0)+\nu(1)+\cdots+\nu(\frac{d}{4})
\ee
where
\be
\nu(k)\eqdef\left\{
\begin{array}{lll}
\lfloor \frac{k}{6} \rfloor & , & \ {\rm if}\ 6|(k-1) \ ;\\
 & & \\
\lfloor \frac{k}{6} \rfloor +1 & , & \ {\rm otherwise}.
\end{array}
\right.
\ee
For $d$ a multiple of $24$, and letting $l=\frac{d}{24}$,
this simplifies to
\be
{\rm dim}\lp
\CC[F]_d^{SL_2}\rp
= 3l^2+3l+1\ .
\ee

\section{The invariant theory computations}
\label{Invthcal}

\subsection{The basic construction}
We now proceed to the definition of the homogeneous version
$\bR$ of the quintic $\caR$: a construction due to Beauville.
In terms of the homogeneous roots
$\la_1,\ldots,\la_5$ in $\PP^1$, write
\be
F(\ux)=(\ux \la_1)\ldots(\ux \la_5)
\ee
and define the five quartics $Q_1,\ldots,Q_5$ by
\be
Q_i(\ux)=(\ux \la_1)\ldots
\widehat{(\ux \la_i)}\ldots
(\ux \la_5)\ .
\ee
Now introduce a new variable $z$ and let
\bea
\bR & \eqdef & \prod\limits_{i=1}^5
\lp
\lp S(Q_i)^3-27T(Q_i)^2  \rp z - S(Q_i)^3
\rp \\
 & = & \sum\limits_{i=0}^5 \cB_i(F) z^{5-i}
\eea
which defines the expressions $\cB_0,\ldots,\cB_5$.

\begin{Lemma}
$\cB_0,\ldots,\cB_5$ are homogeneous polynomial
invariants of $F$, of degree
$24$.
\end{Lemma}

\demo
This is a straightforward application of Proposition \ref{Rootprop}.
To get the degree, one first calculates the weight
by counting the bracket factors $(\la_i \la_j)$ :
\be
\om=5\times 12=60\ ;
\ee
and then uses (\ref{Degweight}) to obtain
\be
d=\frac{2\times 60}{5}=24\ .
\ee

\qed

\subsection{The main results}
The most crucial step in this article is the following
exact determination of the invariants $\cB_0,\ldots,\cB_5$
in terms of $J$, $K$, $L$.

\begin{Proposition}
\label{Keyprop}
\bea
\lefteqn{
\cB_0(F)=\frac{5^{15}}{2^{40}}
\left\{
- 2^{21}. K^3 + 2^{14}.3.K^2 J^2 - 2^7.3.KJ^4 + J^6
\right\}\ ,}  \label{B0formula}\\
\lefteqn{
\cB_1(F)=\frac{5^{16}}{2^{35}.3^3}
\left\{
2^{16}.7.K^3
- 2^{10}.23.K^2 J^2
+ 2^2.71.KJ^4
- J^6
\right\}\ ,} \\
\lefteqn{
\cB_2(F)=\frac{5^{16}}{2^{30}.3^6}
\left\{
2^{11}.5^3.LKJ
- 2^4.5^3.LJ^3
\right.} \nonumber\\
& & \left.
- 2^{15}.3.K^3
+2^7.11.13.K^2 J^2
-3.131.KJ^4
+2.J^6
\right\} \ ,\\
\lefteqn{
\cB_3(F)=\frac{5^{16}}{2^{25}.3^9}
\left\{
-2^{11}.5^4.L^2
- 2^9.3.5^3.LKJ
\right.} \nonumber\\
& & \left.
+ 2.5^3.11.LJ^3
+ 2^9.17.K^3
- 2^2.23.37.K^2J^2
+ 3^5.KJ^4
- 2.J^6
\right\} \ ,\nonumber\\
 & & \\
\lefteqn{
\cB_4(F)=\frac{5^{16}}{2^{22}.3^{12}}
\left\{
-2^5.3^2.5^3.LKJ
- 5^3.29.LJ^3
\right.} \nonumber\\
 & & \left.
-2^7.11.K^3
-7^2.83.K^2 J^2
- 2^2.59.KJ^4
+2^2.J^6
\right\} \ ,\\
\lefteqn{
\cB_5(F)=\frac{5^{15}}{2^{15}.3^{15}}
\left\{
3^3.K^3
-3^3.K^2 J^2
+3^2.KJ^4
- J^6
\right\}\ .}
\eea
\end{Proposition}

\begin{Remark}
It is clear, by construction, that
\be
\cB_0(F)=\prod\limits_{i=1}^5
\lp
2^{-8}
{\rm Disc}(Q_i)
\rp
= 2^{-40}.{\rm Disc}(F)^3
\ee
which can be compared, as a consistency check with (\ref{B0formula})
rewritten as
\be
\cB_0(F)= 2^{-40}
\left[
5^5\lp
J^2-128 K
\rp
\right]^3\ \ .
\ee
Indeed, one can verify, with the help of Maple, that
\be
{\rm Disc}(F)=
5^5\lp
J^2-128 K
\rp\ .
\ee
\end{Remark}

\noindent {\sc Computer-assisted proof of the proposition.}\;
The argument relies on noticing that the construction of $\bR$
is a particular instance of a quartic {\em Tschirnhaus transformation}
of a quintic equation.
Since one already knows that $\cB_0,\ldots,\cB_5$
are homogeneous polynomials of degree $24$ in the coefficients
of the quintic
\be
F(\ux)=
a_0 x_1^5
+ a_1 x_1^4 x_2
+ a_2 x_1^3 x_2^2
+ a_3 x_1^2 x_2^3
+ a_4 x_1 x_2^4
+ a_5 x_2^5\ ;
\ee
one can safely dehomogenize by letting $a_0=1$.
We also dehomogenize with respect to the
variables $x_1$, $x_2$ by letting
$x_1=x$ and $x_2=1$. With a harmless abuse of notation,
the quintic $F$ becomes the monic polynomial
\bea
F(x) & = & x^5+a_1 x^4+ a_2 x^3 +a_3 x^2 + a_4 x + a_5 \\
 & = & (x-\la_1)\ldots(x-\la_5)\ .
\eea
Now the $Q_i$ become
\be
Q_i(x)=\prod\limits_{{j=1}\atop{j\neq i}}^5
(x-\la_i)=\frac{F(x)}{x-\la_i}\ \ .
\ee
In terms of a root $\la$ (or rather a new variable which will
{\em later}
be specialized to such root), the corresponding quartic
is given, after explicitly performing the Euclidean division,
as in (\ref{quartic}) by
\be
Q_\la(x)=
q_0 x^4+ 4 q_1 x^3 + 6 q_2 x^2
+ 4 q_3 x + q_4
\ee
where
\bea
q_0 & = & 1 \ ,\\
q_1 & = & \frac{1}{4} (\la+a_1) \ ,\\
q_2 & = & \frac{1}{6} (\la^2+a_1\la+a_2) \ ,\\
q_3 & = & \frac{1}{4} (\la^3+a_1\la^2+ a_2\la+a_3) \ ,\\
q_4 & = & \la^4+ a_1\la^3 + a_2\la^2 + a_3\la + a_4\ .
\eea
Using the expressions (\ref{Sformula})
and (\ref{Tformula}) for the quartic invariants $S$ and $T$,
one substitutes these values in
\be
\ph(\la)\eqdef
\lp
S(Q_\la)^3-27 T(Q_\la)^2
\rp
z-S(Q_\la)^3\ \ .
\ee
This is, a priori, a polynomial in $\la$ of degree $12$
(i.e., the weight of the isobaric expression
$(S^3-27 T^2)z-S^3$ in the $q$'s).
We now perform the Euclidean division of $\ph(\la)$
by
\be
F(\la) = \la^5+a_1 \la^4+ a_2 \la^3 +a_3 \la^2 + a_4 \la + a_5 \ ,
\ee
and call the remainder $\bar{\ph}(\la)$.
Since the initially generic
$\la$ is going to be specialized to a root of $F$, one will
have
\be
\bR=\prod\limits_{i=1}^5 \ph(\la_i)
=\prod\limits_{i=1}^5 {\bar{\ph}}(\la_i)\ \ .
\ee
By the Poisson product formula this boils down to
\be
\bR={\rm Res}(F,{\bar{\ph}})
\ee
the resultant of two polynomials in $\la$:
$F(\la)$ of degree $5$, and ${\bar{\ph}}(\la)$ of degree $4$.
This is calculated using the Sylvester
determinant formula as in (\ref{Sylvres}).
One obtains $\bR$ as a polynomial in $a_1,\ldots,a_5$, and
$z$. The invariants $\cB_0,\ldots,\cB_5$
are extracted as the coefficients of the powers of $z$.
Now one rehomogenizes by performing the substitutions
\[
(a_1,\ldots,a_5)\rightarrow
\lp
\frac{a_1}{a_0},\ldots,\frac{a_5}{a_0}
\rp\ \ ,
\]
and multiplying by $a_0^{24}$
to get the Cartesian expressions of $\cB_0,\ldots,\cB_5$
as homogeneous polynomials in the
coefficients of the original binary quintic $F$.
Finally one is reduced to a question of linear algebra,
that of decomposing these invariants
in terms of the basis of the degree $24$ component
of the ring $\CC[F]^{SL_2}$ given by the following
monomials in the algebraically independent invariants $J$, $K$, $L$ :
\[
L^2,LKJ,LJ^3,K^3,K^2 J^2, KJ^4,J^6\ .
\]
To make life easier for Maple we did so by first specializing
to the canonical form (\ref{Sylvcan}),
and then solving the linear system in $\CC[u,v,w]$.
The result of these computer calculations is the statement
of the proposition.
\qed

\begin{Remark}
The computationally costly step in this derivation
is the resultant calculation with specialized coefficients
in terms of $a_1,\ldots,a_5$ and $z$.
It took $6$ minutes and $37$ seconds on a
$2\times450$Mhz SUN UltraSparc-II workstation
running Version 9.5  of Maple.
\end{Remark}

Our next computational result is the following.

\begin{Proposition}
\label{Prop48}
The $21$ polynomials $\cB_i^2$, $0\le i\le 5$, and
$\cB_i\cB_j$, $0\le i<j\le 5$, linearly generate
the component of degree $48$ in the
ring of invariants $\CC[F]^{SL_2}$.
\end{Proposition}
\demo
A linear basis of this vector space is given by the $19$
monomials in $J$, $K$, $L$ of that degree.
We simply calculated the $19\times21$ matrix
of coefficients, in this basis, for the $21$ given
polynomials; and we checked, with the help of Maple,
that the matrix has full rank.
\qed

We can now state the main result of this article, which
is the solution to Question \ref{Qinv}.

\begin{Theorem}
\label{Thm48}
For every integer $d>0$ which is a multiple of $d_0=48$,
all invariants of the quintic $F$, of degree $d$, can be
written as polynomials in the invariants $\cB_0,\ldots,\cB_5$.
\end{Theorem}

\demo
Now that Proposition \ref{Prop48} has been established, all one needs to do
is show that for any $d=48k$, where $k\ge 1$ is an integer,
every monomial $L^{\al_1}K^{\al_2}J^{\al_3}$
of degree $d=12\al_1+8\al_2+4\al_3$
can be written as a product of monomials of degree $48$.
This is done by induction on $k$.
For $k=1$, this is a tautology.
Noting that $3\al_1+2\al_2+\al_3=12 k$,
let us perform the Euclidean division of $\al_1$ by $4$,
$\al_2$ by $6$, and $\al_3$ by $12$:
\be
\begin{array}{lll}
\al_1=4\beta_1+\ga_1 & , & 0\le \ga_1\le 3 \ ;\\
\al_2=6\beta_2+\ga_2 & , & 0\le \ga_2\le 5 \ ;\\
\al_3=12\beta_3+\ga_3 & , & 0\le \ga_3\le 11 \ .\\
\end{array}
\label{arithm}
\ee
Clearly, the degree of $L^{4\beta_1}K^{6\beta_2}J^{12\beta_3}$
is a multiple of $48$, and therefore so is that of
$L^{\ga_1}K^{\ga_2}J^{\ga_3}$.
If both triplets $(\beta_1,\beta_2,\beta_3)$
and $(\ga_1,\ga_2,\ga_3)$ are different from $(0,0,0)$,
we are done by induction.

If $(\ga_1,\ga_2,\ga_3)=(0,0,0)$,
then
\be
L^{\al_1}K^{\al_2}J^{\al_3}=
\lp L^4\rp^{\beta_1}
\lp K^6\rp^{\beta_2}
\lp J^{12}\rp^{\beta_3}
\ee
is of the required form.

If $(\beta_1,\beta_2,\beta_3)=(0,0,0)$, then by the
inequalities (\ref{arithm}),
\be
d=12\ga_1+8\ga_2+4\ga_3\le 120\ .
\ee
But $48|d$ and the case $d=48$ has been dealt with; so
we are left with the case where $d=96$.
Since $0\le \ga_1\le 4$, $0\le \ga_2\le 6$ and $3\ga_1+2\ga_2+\ga_3=24$,
one can write
\be
L^{\al_1}K^{\al_2}J^{\al_3}
=L^{\ga_1}K^{\ga_2}J^{\ga_3}
=
\lp
L^{\ga_1} J^{12-3\ga_1}
\rp
\lp
K^{\ga_2} J^{12-2\ga_2}
\rp
\ee
which is the required decomposition.
\qed

\section{The reconstruction problem and the relation to del Pezzo surfaces}
\label{Geom}

The presentation here closely follows, notation included, that of
Beauville~\cite{Beauville1, Beauville2}.

\subsection{The reconstruction problem}

Let $V\eqdef(\PP^1)^5\backslash\De$ where $\De$ is the big diagonal.
One has two commuting actions on $V$ given by that of $SL_2$ and
that of the symmetric group $\fS_5$.
Let $J:V\rightarrow (\PP^1)^5$
be the map which to a quintuple $(\la_1,\ldots,\la_5)$
associates the quintuple $(j_1,\ldots,j_5)$
where $j_l$ is the $j$-invariant, as in (\ref{jinv}), of the
four-point set $\{\la_1,\ldots,
\widehat{j_l},\ldots,j_5\}$.
The map is $\fS_5$-equivariant and factors through
the quotient $P\eqdef V/SL_2$;
i.e., one has a commutative diagram:
\be
\begin{array}{ccc}
P & {\stackrel{J}{\longrightarrow}} & (\PP^1)^5 \\
\downarrow & & \downarrow \\
P/\fS_5 & {\stackrel{\bar{J}}{\longrightarrow}} & {\rm Sym}^5(\PP^1)\ \ .
\end{array}
\ee
The solution to Question \ref{Qreconst} is the following.

\begin{Theorem}
\label{Reconst}
The map $\bar{J}$ is injective.
\end{Theorem}

\demo
Consider two elements $p_1$ and $p_2$ of $P/\fS_5$
which map by $\bar{J}$ to the same element of ${\rm Sym}^5(\PP^1)$.
These correspond to two binary quintics $F_1$ and $F_2$,
defined up to a multiplicative constant.
By hypothesis, the corresponding quintics $\bR$ have the same roots, i.e.,
\be
\forall i,\ 0\le i\le 5,\ \ \ 
\frac{\cB_i(F_1)}{\cB_0(F_1)}
=
\frac{\cB_i(F_2)}{\cB_0(F_2)}\ \ ;
\ee
or what is the same
\be
\forall i,\ 0\le i\le 5,\ \ \ 
\frac{\cB_i(F_1)}{{\rm Disc}(F_1)^3}
=
\frac{\cB_i(F_2)}{{\rm Disc}(F_2)^3}\ \ .
\label{Ratio}
\ee
Now we claim that every expression
$J^\al K^\beta L^\ga H^\de . {\rm Disc}(F)^{-\ep}$ of degree $0$ where
$\al$, $\beta$, $\ga$, $\de$, and $\ep$ are nonnegative
integers, takes the same value for $F_1$ and $F_2$.
Indeed, one has
\be
4\al+8\beta+12\ga+18\de-8\ep=0
\ee
therefore $4|18\de$ so $\de$ is even.
Using the relation (\ref{Relation}) one can get rid of the invariant $H$.
Now
\be
J^\al K^\beta L^\ga . {\rm Disc}(F)^{-\ep}
=
\lp
{\rm Disc}(F)^{\rho} J^\al K^\beta L^\ga 
\rp
{\rm Disc}(F)^{-(\ep+\rho)}
\ee
where $\rho\eqdef
6 \lceil \frac{\ep}{6} \rceil -\ep\ge 0$.
Expressing ${\rm Disc}(F)$ on the left in terms of $J$, and $K$,
one is reduced to the case of an expression
$J^\al K^\beta L^\ga . {\rm Disc}(F)^{-\ep}$ where the degree
of $J^\al K^\beta L^\ga$ is divisible by $48$.
The claim now is a consequence of Theorem \ref{Thm48} and
(\ref{Ratio}).

Now following~\cite[Ch. 5]{Mukai},
$F_1$ and $F_2$ can be seen as elements
of the open affine set
\be
U_{1,5}\eqdef
\left\{
\begin{array}{c}
{\rm binary\ quintics\ without} \\
{\rm repeated\ linear\ factors}
\end{array}
\right\}
\ee
equipped with the natural $GL_2$ action.
It is well-known (see e.g.~\cite[Corollary 5.24]{Mukai})
that $U_{1,5}\rightarrow U_{1,5}/GL_2$ is a good
geometric quotient. The elements of the 
coordinate ring $\CC[F, {\rm Disc}(F)^{-1}]^{GL_2}$
of the latter separate the $GL_2$ orbits.
Using the description of $\CC[F]^{SL_2}$ recalled in Section \ref{Prelim},
these elements are finite linear combinations
of expressions $J^\al K^\beta L^\ga H^\de . {\rm Disc}(F)^{-\ep}$
as above. Now the claim which we have just proved entails:
$F_1$ and $F_2$ are in the same $GL_2$ orbit.
Therefore the corresponding points $p_1$ and $p_2$ in
$P/\fS_5$ are the same.
\qed

\subsection{The relation to del Pezzo surfaces and the Cremona group}

We now come full-circle by explaining how
Theorem \ref{Reconst} provides
an alternate proof of Proposition \ref{Bmapinj}.
The discussion will be quite brief, since much more detail
can be found in~\cite{Beauville1, Beauville2} for the specifics
of the situation, and~\cite{AvritzerM} as well as~\cite[Lecture 22]{Harris}
for the standard prerequisites on quartic del Pezzo surfaces.
Such a surface $S$ is usually seen as the blow up of $\PP^2$
at five points in general position.
The linear system of cubics through these five points embeds
$S$ as a complete intersection of two
quadrics in $\PP^4$.
By choosing an appropriate coordinate system in the latter
one can take these quadrics to be given by
equations $Q_\infty=0$ and $Q_0=0$
where
\be
Q_\infty\eqdef\sum\limits_{i=1}^5 X_i^2\ \ \ {\rm and}\ \ \ 
Q_0\eqdef\sum\limits_{i=1}^5 \la_i X_i^2\ .
\ee
There is a canonical subgroup of automorphisms
of $S$, isomorphic to $(\ZZ/2)^4$, which is the one generated
by the involutions $\si_l$ mapping
$(X_1:\ldots:X_l:\ldots:X_5)$ to $(X_1:\ldots:-X_l:\ldots:X_5)$. 
This descends, via the birational map from $S$ to $\PP^2$
corresponding to the blow up,
to a subgroup $G_S$ isomorphic to $(\ZZ/2)^4$ inside the Cremona group
${\rm Cr}$, or rather to a conjugacy class of such.
This is the construction given by Beauville for the map
in Proposition \ref{Bmapinj}.

Now note that the moduli space of (nonsingular)
quartic del Pezzo surfaces $S$ is the same as that
of binary quintics without repeated linear factors, or more precisely
the space we denoted earlier by $P/\fS_5$.
This correspondence is given by the consideration of the pencil
$Q_\infty \la-Q_0$ which is singular exactly when $\la$
belongs to the five point set $\{\la_1,\ldots,\la_5\}$
(see~\cite{AvritzerM} for
a very thorough treatment).
From the knowledge of the conjugacy class $G_S$ one can
recover the isomorphism class of the
normalized fixed point locus, i.e., the normalization of the union
of the nonrational curves in $\PP^2$ which are fixed by an element
of $G_S$. At the level of the surface $S$, this means that one can
recover the data of the $j$-invariants of the five elliptic curves
obtained as the intersection of $S$ with each of the hyperplanes
$X_l=0$. This is the same as the unordered collection
of the $j_l$'s as in (\ref{jquintic}).
As a result Theorem \ref{Reconst} implies Proposition \ref{Bmapinj}.

\begin{Remark}
We did not try to see if the nice geometric method
used by Beauville in~\cite{Beauville2} could be refined
in order to obtain Theorem \ref{Reconst},
or Theorem \ref{Thm48} (at least with unspecified
$d_0$). This might be an interesting point to elaborate upon
in view of the generalization proposed in Section \ref{TodoAG} below.
\end{Remark}

\section{A shopping list}
\label{Todo}

One of the `raisons d'\^{e}tre' of experiment in natural sciences
is to spur new theoretical investigations.
Accordingly, we would be very happy to see
the experimental mathematical result obtained in this
article initiate some search for theoretical understanding,
however modest.
We can already see different questions arise from this work which
might variously interest the communities of
algebraic geometers,
combinatorial/computational algebraists, and representation theorists.
We will organize these questions accordingly.

\subsection{Algebraic Geometry}
\label{TodoAG}

Very loosely speaking,
our Theorem \ref{Reconst} can be recast
with vast although probably not maximal
generality as the following.
Consider a projective variety $X$ equipped with the action of
a reductive group $G$.
One can try to mimic the construction of Beauville's map $\bar{J}$
and analyse the injectivity of
\be
\mbox{``}\ \ 
{\lp
{\rm Sym}^p X
\rp/G
\rightarrow
\lp
\prod_I
\lp
\lp{\rm Sym}^{|I|} X\rp /G
\rp
\rp/\fS_p}
\ \ \mbox{''}
\ee
where $I$ ranges through the $\binom{p}{q}$ subsets
of cardinality $q$ in $\{1,\ldots,p\}$.
One would have to do some work even in order to give a clean
formulation of the question, in particular with regard
to the analogue of the big diagonal $\De$ one needs to remove
and related stability issues; this is why we put quotes.
In particular, if only for esthetic reasons, one might want
to investigate the case of $SL_{n+1}$ acting on $\PP^n$,
or an invariant theoretic
interrelation of Chow varieties of zero-cycles in $\PP^n$
of different degrees.
A special situation with binary forms, where the precise formulation
of the problem is straightforward is the reconstruction problem
for a binary $p$-ic from the $j$-invariants of its four-root
subsets. It is somewhat the natural one-dimensional projective
analogue of similar questions in distance geometry
and rigidity theory (see e.g.~\cite{BorceaS})
where one tries to determine a Euclidean configuration of points
from mutual distances. Indeed, in the Euclidean situation
{\em one} modulus corresponds to two-point subsets, whereas here
it corresponds to four-point subsets.

\subsection{Combinatorial/Computational Algebra}

A natural problem, under this heading, is to reduce the computations
which we have done
(especially the ones in Proposition \ref{Keyprop}) 
to human proportions.
This might well be needed in order to tackle the next open case of 
reconstructing
the binary
sextic from the quartics it contains.
Indeed, one would have to identify $16$ invariants of degree $60$
one of which is the $6$-th power of the discriminant.
When doing explicit calculations with invariants of
binary forms, one essentially has the following tools
and combinations thereof to choose from:
\begin{enumerate}

\item
Cartesian expressions,

\item
Canonical forms,

\item
Symmetric functions of root differences,

\item
The symbolic method.
\end{enumerate}

In our opinion, the most interesting is the one we did not use
in this article, i.e., the last one.
One would need to {\em invariantively} rephrase
the proof of Proposition \ref{Keyprop}, i.e.,
keeping the $SL_2$-equivariance explicit throughout.
We believe the methods to do that are already available in the
classical literature (see e.g.~\cite{Hermite2, Gordan2, Gordan3, Brioschi})
waiting to be carefully studied anew by computational and combinatorial
algebraists.
Maybe a word of caution for those who would be willing to do so
is in order.
They will find, in addition to the common
vicissitudes of research life, three practical obstacles
specific to this task, and pertaining to

\begin{enumerate}

\item
Physical accessibility of the literature,

\item
Mathematical accessibility of its contents,

\item
Language barrier.

\end{enumerate}

Fortunately, removing the first obstacle is well under way
thanks to the highly commendable efforts of the retrodigitalization
projects throughout the world.
The second one is no obstacle at all, if only psychological.
Indeed, we explained in~\cite{AC1, AC2}
a minimally acrobatic way of making rigorous mathematical
sense out of the symbolic method as used by classical
masters such as P. Gordan.
The third problem is serious and requires a generous,
volunteer-based translation effort following
the example for instance set by Ackerman and Hermann~\cite{Hilbert},
or Cox and Rojas~\cite{Minding}.
Since one should follow one's own advice, let us announce
a forthcoming translation into English by K.~Hoechsmann,
with commentary by the present author, of the
classical masterpiece~\cite{Gordan1}.

\subsection{Representation Theory}

Modern practice in algebraic geometry does not encourage
the writing of equations in coordinates. The successes
obtained in conformity with this ideological choice
can hardly be argued against; the resulting achievements
are among the greatest of $20$-th century mathematics.
However, it is sometimes {\em necessary} to calculate with coordinates,
especially in view of the currently growing importance of
computational/combinatorial algebra.
It is therefore essential, when it is required,
to try to do so wisely; and in this respect,
there is much to be learned from the $19$-th
century mathematicians.
To continue on what we said in the previous subsection, one
has to realize that from the mere use of Cartesian expressions
one is automatically breaking $SL$-invariance and, prehaps unwittingly,
doing toric geometry.
More appropriate to calculations in the realm of
projective geometry is the symbolic method which
{\em explictly} preserves the $SL$-equivariance.
Concerning the latter, there are a few questions
arising as to what is the representation theoretic
interpretation of our Theorem \ref{Thm48}.

Let $S_d(\cdot)$ denote the $d$-th symmetric power
of an $SL_2$ representation; if no argument is indicated
it means that of the defining vector space $\CC^2$ which
is also identified with its dual.
It is not hard to rephrase our Theorem \ref{Thm48}
as the surjectivity of an $SL_2$-equivariant map
\[
S_{2k}\lp
S_5\lp
\left[
S_6(S_4)
\right]^{SL_2}
\rp
\rp
\rightarrow
\lp
S_{48k}(S_5)
\rp^{SL_2}
\]
One can then ask if one could remove the restriction
to $SL_2$-invariant subspaces.
Indeed the construction of the invariants $\cB_0,\ldots,\cB_5$
is susceptible of many variations and twists.
For instance, one can do it not only for
invariants but also for {\em covariants}
since Proposition \ref{Rootprop} works
equally well for them.
This ties in with one of the main themes of the
article~\cite{AC1}, as well as a rather mysterious
`devissage' property
of classical invariants alluded to in~\cite[pp.~114-118]{Sylvester}.
As an exercise we leave to the reader, and as an illustration
of the point we are making,
one can do the following construction.
Take the covariant $C(Q;\ux)=(Q,(Q,Q)_2)_1$ of the
quartic $Q$; and similarly to the description of $\bR$,
consider the expression
\[
\prod\limits_{i=1}^5
C(Q_i; \xi_i)
\]
which reintroduces the missing root by specializing $\ux$.
Now one can check that this is a nonzero numerical
multiple of the invariant $H$.
This gives a somewhat less `out-of-the-blue'
derivation for the root-difference
expression found by Hermite for his own invariant $H$
(see~\cite{Hermite3} and e.g.~\cite{Kraft}
for related recent work).

\vskip 1cm

{\sc Acknowledgements.} {\small We are very grateful to
David Brydges and Joel Feldman for their invitation to the
University of British Columbia.
We thank Jaydeep Chipalkatti for sharing some of his Maple routines and
for useful discussions. Discussions with Zinovy Reichstein were also
very useful. We are very grateful to Arnaud Beauville
for helping us reach, through email correspondence, the correct
formulation of Question \ref{Qinv}.
We were impressed by how the Maple software handled the resultant
calculation in Proposition \ref{Keyprop}.
The following electronic libraries have
been useful in accessing classical
references:

\begin{itemize}

\item Gallica, Biblioth\`eque Nationale de France ({\bf GA}),

\item The G\"ottinger DigitalisierungsZentrum ({\bf GDZ}),

\item JSTOR ({\bf JS}),

\item The University of Michigan Historical Mathematics Collection ({\bf UM}).

\end{itemize}

\bibliographystyle{plain}

\begin{thebibliography}{10}

\bibitem{AC1}
A.~Abdesselam and J.~Chipalkatti.
\newblock Brill-Gordan loci, transvectants and an analogue
of the Foulkes conjecture .
\newblock preprint math.AG/0411110, 2004.

\bibitem{AC2}
A.~Abdesselam and J.~Chipalkatti.
\newblock The bipartite Brill-Gordan locus and angular momentum.
\newblock preprint math.AG/0502542, 2005.

\bibitem{AvritzerM}
D.~Avritzer and R.~Miranda.
\newblock Stability of pencils of quadrics in $\PP\sp 4$.
\newblock {\em Bol. Soc. Mat. Mexicana (3)}, 
vol.~5, No.~2, pp.~281--300, 1999.

\bibitem{Beauville1}
A.~Beauville.
\newblock $p$-elementary subgroups of the Cremona group.
\newblock preprint math.AG/0502123 v1, 2005.

\bibitem{Beauville2}
A.~Beauville.
\newblock $p$-elementary subgroups of the Cremona group.
\newblock preprint http://math1.unice.fr/~beauvill/pubs/pp.pdf, 2005.

\bibitem{BorceaS}
C.~Borcea and I.~Streinu.
\newblock The number of embeddings of minimally rigid graphs.
\newblock {\em Discrete Comput. Geom.}, vol.~31, No.~2, pp.~287--303, 2004.

\bibitem{Brioschi}
F.~Brioschi.
\newblock Sur la transformation des \'equations alg\'ebriques.
\newblock {\em C. R. Acad. Sci.}, vol.~124 , pp.~661--665, 1897
({\bf GA}).

\bibitem{Elliott}
E.~B.~Elliott.
\newblock {\em An Introduction to the Algebra of Quantics}.
\newblock 2nd ed., Oxford University Press, Oxford, 1913.

\bibitem{Gordan1}
P.~Gordan.
\newblock Beweis, dass jede Covariante und Invariante einer
bin\"aren Form eine ganze Function mit
numerischen Coefficienten einer endlichen Anzahl
solcher Formen ist.
\newblock {\em J.~Reine~Angew.~Math.}, vol.~69, pp.~323--354, 1868
({\bf GDZ}).

\bibitem{Gordan2}
P.~Gordan.
\newblock Ueber die Invarianten bin\"{a}rer Formen bei
h\"{o}heren Transformationen.
\newblock {\em J.~Reine~Angew.~Math.}, vol.~71, pp.~164--194, 1870
({\bf GDZ}).

\bibitem{Gordan3}
P.~Gordan.
\newblock Ueber die Bildung der Resultante zweier Gleichungen.
\newblock {\em Math. Ann.}, vol.~3 , pp.~ 355--414, 1871
({\bf GDZ}).

\bibitem{Harris}
J.~Harris.
\newblock {\em Algebraic Geometry, A First Course}.
\newblock Graduate Texts in Mathematics. Springer--Verlag, \, New York, 1992.

\bibitem{Hermite1}
C.~Hermite.
\newblock
Sur le th\'eorie des fonctions homog\`enes \`a deux ind\'etermin\'ees.
\newblock {\em Cambridge and Dublin Math. J.}, vol.~9, pp.~172--217,
1854 ({\bf UM}).

\bibitem{Hermite2}
C.~Hermite.
\newblock Sur quelques th\'eor\`emes d'alg\`ebre et la
r\'esolution de l'\'equation du quatri\`eme degr\'e.
\newblock {\em C. R. Acad. Sci.}, vol.~46 , pp.~961--967, 1856
({\bf GA, UM}).

\bibitem{Hermite3}
C.~Hermite.
\newblock Sur l'invariant du $18^{\rm e}$ ordre des formes
du cinqui\`eme degr\'e et sur le role qu'il joue dans la
r\'esolution de l'\'equation du cinqui\`eme degr\'e,
extrait de deux lettres de M. Hermite \`a l'\'editeur.
\newblock {\em J.~Reine~Angew.~Math.}, vol.~59, pp.~304--305, 1861
({\bf GDZ, UM}).

\bibitem{Hilbert}
D.~Hilbert.
\newblock {\em Hilbert's Invariant Theory Papers}. 
\newblock Translated from the German by Michael Ackerman.
With comments by Robert Hermann.
Lie Groups: History, Frontiers and Applications, VIII. 
Math. Sci. Press, Brookline, Mass., 1978. 

\bibitem{Iskovskikh}
V.~A.~Iskovskikh.
\newblock Rational surfaces with a pencil of rational curves and
with positive square of the canonical class.
\newblock {\em Math. USSR Sbornik}, vol.~12, pp.~93--117, 1970.

\bibitem{Kempe}
A.~B.~Kempe.
\newblock On the application of Clifford's graphs to ordinary
binary quantics.
\newblock {\em Proc.~London Math.~Soc.}, vol.~17, pp.~107--121, 1885.

\bibitem{Kraft}
H.~Kraft.
\newblock A result of Hermite and equations of degree 5 and 6.
\newblock {\em J. Algebra}, in press, 2005. 

\bibitem{Kung}
J.P.S. Kung.
\newblock Canonical forms of binary forms: variations on
a theme of Sylvester. 
\newblock {\em Invariant theory and tableaux (Minneapolis, MN, 1988)}, 
pp.~46--58, IMA Vol. Math. Appl., 19, 
Springer Verlag, New York, 1990. 

\bibitem{KungR}
J.P.S. Kung and G.-C. Rota.
\newblock The invariant theory of binary forms.
\newblock {\em Bulletin of the A.M.S.}, vol. 10, \,No. 1,\,pp. 27--85, 1984.

\bibitem{Minding}
F.~Minding.
\newblock On the determination of the degree of an equation
obtained by elimination. 
\newblock Translated from the German and with
a commentary by D. Cox and J. M. Rojas. Contemp. Math., 334, 
{\em Topics in Algebraic Geometry and Geometric Modeling}, pp.~351--362, 
Amer. Math. Soc., Providence, RI, 2003. 

\bibitem{Mukai}
S.~Mukai.
\newblock {\em An Introduction to Invariants and Moduli}.
\newblock Translated from the 1998 and 2000 Japanese 
editions by W. M. Oxbury. Cambridge Studies in Advanced Mathematics, 81. 
Cambridge University Press, Cambridge, 2003.

\bibitem{Salmon}
G.~Salmon.
\newblock {\em Higher Algebra} (5th ed.), 1885.
\newblock Reprinted by Chelsea Publishing Co., New York, 1964.

\bibitem{Sylvester}
J.~J.~Sylvester.
\newblock On an application of the new atomic theory to the graphical
representation of the invariants and covariants of binary quantics,
with three appendices.
\newblock {\em  Amer.~J.~Math.}, vol.~1, pp.~64--125, 1878 ({\bf JS}).

\end{thebibliography}

\bigskip
\centerline{------------------------------}

\vspace{3cm}

\parbox{6cm}{\small 
{\sc Abdelmalek Abdesselam} \\ 
Department of Mathematics\\
University of British Columbia\\
1984 Mathematics Road \\
Vancouver, BC V6T 1Z2 \\ Canada. \\ 
{\tt abdessel@math.ubc.ca}} 
\hfill 
\parbox{5cm}{\small 
LAGA, Institut Galil\'ee \\ CNRS UMR 7539\\
Universit{\'e} Paris XIII\\
99 Avenue J.B. Cl{\'e}ment\\
F93430 Villetaneuse \\ France. \\
{\tt abdessel@math.univ-paris13.fr}}

\end{document}